*Research Article*

# Algorithm for Solutions of Nonlinear Equations of Strongly Monotone Type and Applications to Convex Minimization and Variational Inequality Problems

Mathew O. Aibinu ,[1] Surendra C. Thakur,[2] and Sibusiso Moyo[3]

[1]*Institute for Systems Science & KZN CoLab, Durban University of Technology, Durban 4000, South Africa*
[2]*KZN CoLab, Durban University of Technology, Durban 4000, South Africa*
[3]*Institute for Systems Science & Office of the DVC Research, Innovation & Engagement Milena Court, Durban University of Technology, Durban 4000, South Africa*

Correspondence should be addressed to Mathew O. Aibinu; moaibinu@yahoo.com

Real-life problems are governed by equations which are nonlinear in nature. Nonlinear equations occur in modeling problems, such as minimizing costs in industries and minimizing risks in businesses. A technique which does not involve the assumption of existence of a real constant whose calculation is unclear is used to obtain a strong convergence result for nonlinear equations of $(p,\eta)$-*strongly monotone* type, where $\eta > 0, p > 1$. An example is presented for the nonlinear equations of $(p,\eta)$-*strongly monotone* type. As a consequence of the main result, the solutions of convex minimization and variational inequality problems are obtained. This solution has applications in other fields such as engineering, physics, biology, chemistry, economics, and game theory.

## 1. Introduction

Let $\nu : [0,\infty) \longrightarrow [0,\infty)$ be a continuous, strictly increasing function such that $\nu(t) \longrightarrow \infty$ as $t \to \infty$ and $\nu(0) = 0$ for any $t \in [0,\infty)$. Such a function $\nu$ is called a *gauge function*. Let $B$ be a Banach space, and $B^*$ denotes it is dual. A *duality mapping* associated with the gauge function $\nu$ is a mapping $J_\nu^B : B \longrightarrow 2^{B^*}$ defined by

$$J_\nu^B(x) = \{f \in B^* : \langle x, f \rangle = \|x\|\nu(\|x\|), \|f\| = \nu(\|x\|)\}, \quad (1)$$

where $\langle ., . \rangle$ denotes the duality pairing. For $p > 1$, let $\nu(t) = t^{p-1}$ be a gauge function. $J^B : B \longrightarrow 2^{B^*}$ is called a *generalized duality mapping* from $B$ into $2^{B^*}$ and is given by

$$J_p^B(x) = \left\{ f \in B^* : \langle x, f \rangle = \|x\|^p, \|f\| = \|x\|^{p-1} \right\}. \quad (2)$$

For $p = 2$, the mapping $J_2^B$ is called the *normalized duality mapping*. In a Hilbert space, the normalized duality mapping is the identity map. For $x, y \in B$ and $\eta > 0$, a mapping $T : B \longrightarrow B^*$ is said to be

(1) *monotone* if $\langle x - y, Tx - Ty \rangle \geq 0$

(2) *strongly monotone* (see, e.g., Alber and Ryazantseva [1], p. 25), if $\langle x - y, Tx - Ty \rangle \geq \eta\|x - y\|^2$

(3) $\eta$-*strongly pseudomonotone* if $\langle x - y, Ty \rangle \geq 0 \implies \langle x - y, Tx \rangle \geq \eta\|x - y\|^2$

(4) $(p, \eta)$-*strongly monotone* if $\langle x - y, Tx - Ty \rangle \geq \eta\|x - y\|^p$ (see, e.g., Chidume and Djitté [2], Chidume and Shehu [3], and Aibinu and Mewomo [4, 5]).

*Remark 1.* According to the definition of Chidume and Djitté [2] and Chidume and Shehu [3], a strongly monotone mapping is referred to as $(2, \eta)$-strongly monotone mapping.

A monotone mapping $T$ is called *maximal monotone* if it is a monotone and its graph is not properly contained in the graph of any other monotone mapping. As a result of Rockafellar [6], $T$ is said to be a maximum monotone if it is a monotone and the range of $(J^B + tT)$ is all of $B^*$ for some $t > 0$. The set of zeros of a maximum monotone mapping $T$, $T^{-1}(0) := \{x \in B : Tx = 0\}$ is closed and convex. A function $\varphi : B \longrightarrow (-\infty, +\infty]$ is said to be *proper* if the set $\{x \in \mathbb{R} : F(x) \in \mathbb{R}\}$ is nonempty. A proper function $\varphi : B \longrightarrow (-\infty, +\infty]$ is said to be *convex* if for all $x, y \in B$ and $\tau \in [0, 1]$, we have

$$\varphi(\tau x + (1 - \tau)y) \leq \tau\varphi(x) + (1 - \tau)\varphi(y). \quad (3)$$

If the set of $\{x \in \mathbb{R} : \varphi(x) \leq r\}$ is closed in $B$ for all $r \in \mathbb{R}, \varphi$ is said to be *lower semicontinuous*. For a proper lower semicontinuous function $\varphi : B \longrightarrow (-\infty, +\infty]$, the *subdifferential mapping* $\partial\varphi : B \longrightarrow 2^{B^*}$, defined by

$$\partial\varphi(x) = \{x^* \in B^* : \varphi(y) - \varphi(x) \geq \langle y - x, x^* \rangle \forall y \in B\}, \quad (4)$$

is a maximal monotone (Rockafellar [7]). Consider a problem of finding a solution of the equation $Tu = 0$, where $T$ is a maximal monotone mapping. Such a problem is associated with the *convex minimization problem*. Indeed, for a proper lower semicontinuous convex function $\varphi : B \longrightarrow (-\infty, +\infty]$, solving the equation $Tu = 0$ is equivalent to finding $\varphi(u) = \min_{x \in B} \varphi(x)$ by setting $\partial\varphi \equiv T$.

For a reflexive smooth strictly convex space $B$, we let $T$ be a mapping such that the range of $(J_p^B + tT)$ is all of $B^*$ for some $t > 0$ and let $x \in B$ be fixed. Then, for every $t > 0$, there corresponds a unique element $x_t \in D(T)$ such that

$$J_p^B x = J_p^B x_t + tTx_t. \quad (5)$$

Therefore, the *resolvent* of $T$ is defined by $J_t^T x = x_t$. In other words, $J_t^T = (J_p^B + tT)^{-1} J_p^B$ and $T^{-1}0 = F(J_t^T)$ for all $t > 0$, where $F(J_t^T)$ denotes the set of all fixed points of $J_t^T$. The resolvent $J_t^T$ is a single-valued mapping from $B$ into $D(T)$ (Kohsaka and Takahashi [8]). $J_t^T$ is nonexpansive if $E$ is a Hilbert space (Takahashi [9]). Some existing results proved a strong convergence theorem for nonlinear equations of the monotone type, with the assumption of existence of a real constant whose calculation is unclear (see, e.g., Aibinu and Mewomo [4], Chidume et al. [10], and Diop et al. [11]). Monotone-type mappings occur in many functional equations, and the research on monotone type mappings has recently attracted much attention (see, e.g., Shehu [12, 13], Chidume et al. [14], Djitte et al. [15], Tang [16], Uddin et al. [17], Chidume and Idu [18], and Aibinu and Mewomo [19]).

In this paper, we consider nonlinear equations of $(p, \eta)$-strongly monotone type, $p > 1$ and $\eta \in (1, \infty)$. This is a wider class than the class of strongly monotone mappings. An example is presented for nonlinear equations of $(p, \eta)$-*strongly monotone* type. Under suitable conditions which do not involve the assumption of existence of a real constant whose calculation is unclear, a sequence of iteration is shown to converge strongly to the zero of a nonlinear equation of $(p, \eta)$-strongly monotone type. As a consequence of the main result, the solution of convex minimization and variational inequality problems is obtained, which has applications in several fields such as economics, game theory, and the sciences.

## 2. Preliminaries

Let $B$ be a real Banach space and $S := \{x \in B : \|x\| = 1\}$. $B$ is said to have a *Gateaux differentiable norm* if the limit

$$\lim_{t \to 0} \frac{\|x + ty\| - \|x\|}{t}, \quad (6)$$

exists for each $x, y \in S$. A Banach space $B$ is said to be *smooth* if for every $x \neq 0$ in $B$, there is a unique $x^* \in B^*$ such that $\|x^*\| = 1$ and $\langle x, x^* \rangle = \|x\|$, where $B^*$ denotes the dual of $B$. $B$ is said to be *uniformly smooth* if it is smooth and the limit (6) is attained uniformly for each $x, y \in S$. The *modulus of convexity* of a Banach space $B$, $\delta_B : (0, 2] \longrightarrow [0, 1]$ is defined by

$$\delta_B(\varepsilon) = \inf\left\{1 - \frac{\|x + y\|}{2} : \|x\| = \|y\| = 1, \|x - y\| > \varepsilon\right\}. \quad (7)$$

$B$ is *uniformly convex* if and only if $\delta_B(\varepsilon) > 0$ for every $\varepsilon \in (0, 2]$. A normed linear space $B$ is said to be *strictly convex* if

$$\|x\| = \|y\| = 1, x \neq y \Longrightarrow \frac{\|x + y\|}{2} < 1. \quad (8)$$

It is well known that a space $B$ is uniformly smooth if and only if $B^*$ is uniformly convex.

A mapping $T : B \longrightarrow B^*$ is locally bounded at $v \in D$, if there exist $r_v > 0$ and $m > 0$ such that

$$\|Tx\| \leq m, \forall x \in D_{r_v}(v). \quad (9)$$

In particular, $\|Tv\| \leq m$. Therefore, $\langle v, Tv \rangle \leq m\|v\|$. Let $X$ and $Y$ be real Banach spaces and let $T : X \longrightarrow Y$ be a mapping. $T$ is *uniformly continuous* if for each $\varepsilon > 0$, there exists $\delta > 0$ such that

$$\|Tx - Ty\| < \varepsilon \quad \forall x, y \in X \text{ with } \|x - y\| < \delta. \quad (10)$$

Let $\psi(t)$ be a function on the set $\mathbb{R}^+$ of nonnegative real numbers such that

(i) $\psi$ is nondecreasing and continuous

(ii) $\psi(t) = 0$ if and only if $t = 0$

$T$ is said to be uniformly continuous if it admits the modulus of continuity $\psi$ such that

$$\|T(x) - T(y)\| \leq \psi(\|x - y\|) \quad \forall x, y \in X. \quad (11)$$

The modulus of continuity $\psi$ has some useful properties (for instance, see Altomare and Campiti [20], pp. 266–269; Forster [21] and references therein).

2.1. *Properties of Modulus of Continuity.* Let $X$ and $Y$ be real Banach spaces and let $T : X \longrightarrow Y$ be a map which admits the modulus of continuity $\psi$.

(a) *Modulus of continuity is subadditive*: for all real numbers $t_1 \geq 0, t_2 \geq 0$, we have

$$\psi(t_1 + t_2) \leq \psi(t_1) + \psi(t_2) \tag{12}$$

(b) *Modulus of continuity is monotonically increasing*: if $0 \leq t_1 \leq t_2$ holds for some real numbers $t_1, t_2$, then

$$0 \leq \psi(t_1) \leq \psi(t_2) \tag{13}$$

(c) *Modulus of continuity is continuous*: the modulus of continuity $\psi : \mathbb{R}^+ \longrightarrow \mathbb{R}^+$ is continuous on the set positive real numbers; in particular, the limit of $\psi$ at 0 from above is

$$\lim_{t \to 0} \psi(t) = 0 \tag{14}$$

Let $C$ be a nonempty subset of a Banach space $B$ and $T$ be a mapping from $C$ into itself.

(i) $T$ is nonexpansive provided $\|Tx - Ty\| \leq \|x - y\|$ for all $x, y \in C$

(ii) $T$ is firmly nonexpansive type (see, e.g., [22]) if $\langle Tx - Ty, j_p^B Tx - j_p^B Ty \rangle \leq \langle Tx - Ty, j_p^B x - j_p^B y \rangle$ for all $x, y \in C$ and $j_p^B \in J_p^B$

The following results about the generalized duality mappings are well known which are established in, e.g., Alber and Ryazantseva [1] (p. 36), Cioranescu [23] (pp. 25–77), Xu and Roach [24], and Zalinescu [25]. Let $B$ be a Banach space.

(i) $B$ is smooth if and only if $J_p^B$ is single-valued

(ii) If $B$ is reflexive, then $J_p^B$ is onto

(iii) If $B$ has uniform Gateaux differentiable norm, then $J_p^B$ is norm-to-weak∗ uniformly continuous on bounded sets

(iv) $B$ is uniformly smooth if and only if $J_p^B$ is single-valued and uniformly continuous on any bounded subset of $B$

(v) If $B$ is strictly convex, then $J_p^B$ is one-to-one, that is, $\forall x, y \in B, x \neq y \Longrightarrow J_p^B(x) \cap J_p^B(y) = \emptyset$

(vi) If $B$ and $B^*$ are strictly convex and reflexive, then $J_p^{B^*}$ is the generalized duality mapping from $B^*$ to $B$, and $J_p^{B^*}$ is the inverse of $J_p^B$

(vii) If $B$ is uniformly smooth and uniformly convex, the generalized duality mapping $J_p^{B^*}$ is uniformly continuous on any bounded subset of $B^*$

(viii) If $B$ and $B^*$ are strictly convex and reflexive, for all $x \in B$ and $f \in B^*$, the equalities $J_p^B J_p^{B^*} f = f$ and $J_p^{B^*} J_p^B x = x$ hold

*Definition 2.* Alber [26] introduced the functions $\phi : B \times B \to \mathbb{R}$, defined by

$$\phi(x, y) = \|x\|^2 - 2\langle x, J_2^B y \rangle + \|y\|^2, \quad \text{for all } x, y \in B, \tag{15}$$

where $J_2^B$ is the normalized duality mapping from $B$ to $B^*$. Let $B$ be a smooth real Banach space and $p, q > 1$ with $1/p + 1/q = 1$. Aibinu and Mewomo [4] introduced the functions $\phi_p : B \times B \longrightarrow \mathbb{R}$, defined by

$$\phi_p(x, y) = \frac{p}{q}\|x\|^q - p\langle x, J_p^B y \rangle + \|y\|^p, \quad \text{for all } x, y \in B, \tag{16}$$

and $V_p : B \times B^* \to \mathbb{R}$, defined as

$$V_p(x, x^*) = \frac{p}{q}\|x\|^q - p\langle x, x^* \rangle + \|x^*\|^p \quad \forall x \in B, x^* \in B^*, \tag{17}$$

where $J_p^B$ is the generalized duality mapping from $B$ to $B^*$.

*Remark 3.* These remarks follow from Definition 2:

(i) For $p = 2$, $\phi_2(x, y) = \phi(x, y)$, which is the definition of Alber [26]. It is easy to see from the definition of the function $\phi$ that

$$(\|x\| - \|y\|)^2 \leq \phi(x, y) \leq (\|x\| + \|y\|)^2 \quad \text{for all } x, y \in B. \tag{18}$$

Indeed,

$$(\|x\| - \|y\|)^2 = \|x\|^2 - 2\|x\|\|y\| + \|y\|^2 \leq \|x\|^2 - 2\langle x, J_p^B y \rangle + \|y\|^2$$
$$= \phi(x, y) \leq \|x\|^2 + 2\|x\|\|y\| + \|y\|^2 = (\|x\| + \|y\|)^2. \tag{19}$$

By similar analysis, it can verified that for each $p \geq 2$,

$$(\|x\|-\|y\|)^p \leq \phi_p(x,y) \leq (\|x\|+\|y\|)^p \quad \text{for all } x, y \in B. \quad (20)$$

(ii) It is obvious that

$$V_p(x, x^*) = \phi_p\left(x, J_p^{B^*} x^*\right) \quad \forall x \in B, x^* \in B^*. \quad (21)$$

Let $B$ be a topological real vector space and $T$ a multivalued mapping from $B$ into $2^{B^*}$. Cauchy-Schwartz's inequality is given by

$$|\langle x, y^* \rangle| \leq \langle x, x^* \rangle^{1/2} \langle y, y^* \rangle^{1/2}, \quad (22)$$

for any $x$ and $y$ in $B$ and any choice of $x^* \in Tx$ and $y^* \in Ty$ (Zarantonello [27]).

In the sequel, we shall need the lemmas whose proofs have been established (see, e.g., Alber [26] and Aibinu and Mewomo [4]).

**Lemma 4.** *Let $B$ be a strictly convex and uniformly smooth real Banach space and $p > 1$. Then,*

$$V_p(x, x^*) + p\left\langle J_p^{B^*} x^* - x, y^* \right\rangle \leq V_p(x, x^* + y^*), \quad (23)$$

*for all $x \in B$ and $x^*, y^* \in B^*$.*

**Lemma 5.** *Let $B$ be a smooth uniformly convex real Banach space and $p > 1$ be an arbitrarily real number. For $d > 0$, let $B_d(0) \coloneqq \{x \in B : \|x\| \leq d\}$. Then, for arbitrary $x, y \in B_d(0)$,*

$$\|x-y\|^p \geq \phi_p(x,y) - \frac{p}{q}\|x\|^q, \quad \text{where } \frac{1}{p} + \frac{1}{q} = 1. \quad (24)$$

**Lemma 6.** *Let $B$ be a reflexive strictly convex and smooth real Banach space and $p > 1$. Then,*

$$\begin{aligned}\phi_p(y,x) - \phi_p(y,z) &\geq p\left\langle z-y, J_p^B x - J_p^B z \right\rangle \\ &= p\left\langle y-z, J_p^B z - J_p^B x \right\rangle \quad \text{for all } x, y, z \in B.\end{aligned} \quad (25)$$

**Lemma 7.** *Let $B$ be a real uniformly convex Banach space. For arbitrary $r > 0$, let $B_r(0) \coloneqq \{x \in B : \|x\| \leq r\}$. Then, there exists a continuous strictly increasing convex function*

$$g : [0, \infty) \longrightarrow [0, \infty), g(0) = 0, \quad (26)$$

*such that for every $x, y \in B_r(0), j_p^B(x) \in J_p^B(x), j_p^B(y) \in J_p^B(y)$, we have $\langle x-y, j_p^B(x) - j_p^B(y) \rangle \geq g(\|x-y\|)$ (see Xu [28]).*

**Lemma 8.** *Let $\{a_n\}$ be a sequence of nonnegative real numbers satisfying the following relations:*

$$a_{n+1} \leq (1 - \alpha_n)a_n + \alpha_n \sigma_n + \gamma_n, n \in \mathbb{N}, \quad (27)$$

*where*

(i) $\{\alpha\}_n \subset (0,1), \sum_{n=1}^{\infty} \alpha_n = \infty$

(ii) $\limsup\{\sigma\}_n \leq 0$

(iii) $\gamma_n \geq 0, \sum_{n=1}^{\infty} \gamma_n < \infty$

*Then, $a_n \longrightarrow 0$ as $n \longrightarrow \infty$ (see Xu [29]).*

**Lemma 9.** *Let $B$ be a smooth uniformly convex real Banach space and let $\{x_n\}$ and $\{y_n\}$ be two sequences from $B$. If either $\{x_n\}$ or $\{y_n\}$ is bounded and $\phi(x_n, y_n) \longrightarrow 0$ as $n \longrightarrow \infty$, then $\|x_n - y_n\| \longrightarrow 0$ as $n \longrightarrow \infty$ (see Kamimura and Takahashi [30]).*

**Lemma 10.** *A monotone map $T : B \longrightarrow B^*$ is locally bounded at the interior points of its domain (see, e.g., Rockafellar [31] and Pascali and Sburlan [32]).*

**Lemma 11.** *If a functional $\phi$ on the open convex set $M \subset \text{dom } \phi$ has a subdifferential, then $\phi$ is convex and lower semicontinuous on the set (see Alber and Ryazantseva [1], p. 17).*

**Lemma 12.** *Let $X$ and $Y$ be real normed linear spaces and let $T : X \longrightarrow Y$ be a uniformly continuous map. For arbitrary $r > 0$ and fixed $x^* \in X$, let*

$$B_X(x^*, r): \{x \in X : \|x - x^*\|_X \leq r\}. \quad (28)$$

*Then, $T(B(x^*, r))$ is bounded (see, e.g., Chidume and Djitte [33]).*

## 3. Main Results

**Theorem 13.** *Let $B$ be a uniformly smooth and uniformly convex real Banach space. Let $p > 1, \eta \in (1, \infty)$; suppose $T : B \longrightarrow B^*$ is a continuous $(p, \eta)$-strongly monotone mapping such that the range of $(J_p^B + tT)$ is all of $B^*$ for all $t > 0$ and $T^{-1}(0) \neq \emptyset$. Let $\{\lambda_n\}_{n=1}^{\infty} \subset (0,1)$ and $\{\theta_n\}_{n=1}^{\infty}$ in $(0, 1/2)$ be real sequences such that*

(i) $\lim_{n \longrightarrow \infty} \theta_n = 0$ and $\{\theta_n\}_{n=1}^{\infty}$ is decreasing

(ii) $\sum_{n=1}^{\infty} \lambda_n \theta_n = \infty$

(iii) $\lim_{n \longrightarrow \infty} ((\theta_{n-1}/\theta_n) - 1)/\lambda_n \theta_n = 0, \sum_{n=1}^{\infty} \lambda_n < \infty \forall n \in \mathbb{N}$

*For arbitrary $x_1 \in B$, define $\{x_n\}_{n=1}^{\infty}$ iteratively by:*

$$x_{n+1} = J_p^{B^*}\left(J_p^B x_n - \lambda_n \left(Tx_n + \theta_n \left(J_p^B x_n - J_p^B x_1\right)\right)\right), \quad n \in \mathbb{N}, \quad (29)$$

where $J_p^B$ is the generalized duality mapping from $B$ into $B^*$. Then, the sequence $\{x_n\}_{n=1}^{\infty}$ converges strongly to the solution of $Tx = 0$.

*Proof.* Observe that there is no need for constructing a convergence sequence if $x = 0$ because it is a zero of $T$ (since $T$ is strongly monotone, which is one to one). Consequently, we are looking for a unique nonzero solution of $Tx = 0$. The proof is divided into two parts.

*Part 1*: the sequence $\{x_n\}_{n=1}^{\infty}$ is shown to be bounded.

Let $q > 1$ with $1/p + 1/q = 1$ and $x \in B$ be a solution of the equation $Tx = 0$. It suffices to show that $\phi_p(x, x_n) \leq r, \forall n \in \mathbb{N}$. The induction method will be adopted. Let $r > 0$ be sufficiently large such that

$$r \geq \max\left\{\phi_p(x, x_1), 4M_0 M, \frac{4p}{q}\|x\|^q\right\}, \quad (30)$$

where $M_0 > 0$ and $M > 0$ are arbitrary but fixed. For $n = 1$, by construction, we have that $\phi_p(x, x_1) \leq r$ for real $p > 1$. Assume that $\phi_p(x, x_n) \leq r$ for some $n \geq 1$. From inequality (20), we have $\|x_n\| \leq r^{1/p} + \|x\|$. Let $D := \{z \in B : \phi_p(x, z) \leq r\}$. Next is to show that $\phi_p(x, x_{n+1}) \leq r$. It is known that $T$ is locally bounded (Lemma 10) and $J_p^B$ is uniformly continuous on bounded subsets of $B$. Define

$$M_0 := \sup\left\{\|Tx_n + \theta_n\left(J_p^B x_n - J_p^B x_1\right)\| : \theta_n \in \left(0, \frac{1}{2}\right), x_n \in D\right\} + 1. \quad (31)$$

Let $\psi$ denotes the modulus of continuity of $J_p^{B^*}$. Then,

$$\begin{aligned}\|x_n - x_{n+1}\| &= \|x_n - J_p^{B^*}\left(J_p^B x_n - \lambda_n\left(Tx_n + \theta_n\left(J_p^B x_n - J_p^B x_1\right)\right)\right)\| \\ &= \|J_p^{B^*}\left(J_p^B x_n\right) - J_p^{B^*}\left(J_p^B x_n - \lambda_n\left(Tx_n + \theta_n\left(J_p^B x_n - J_p^B x_1\right)\right)\right)\| \\ &\leq \psi\left(|\lambda_n|\|Tx_n + \theta_n\left(J_p^B x_n - J_p^B x_1\right)\|\right) \\ &\leq \psi(|\lambda_n|M_0) \leq \psi(\sup\{|\lambda_n|M_0 : \lambda_n \in (0,1)\}).\end{aligned} \quad (32)$$

Since $T$ is locally bounded and the duality mapping $J_p^B$ is uniformly continuous on bounded subsets of $B$, the $\sup\{|\lambda_n|M_0\}$ exists and it is a real number different from infinity. Choose $M := \psi(\sup\{|\lambda_n|M_0\})$. Applying Lemma 4 with $y^* := \lambda_n(Tx_n + \theta_n(J_p^B x_n - J_p^B x_1))$ and by using the definition of $x_{n+1}$, we compute as follows:

$$\begin{aligned}\phi_p(x, x_{n+1}) &= \phi_p\left(x, J^{B^*}\left(J_p^B x_n - \lambda_n\left(Tx_n + \theta_n\left(J_p^B x_n - J_p^B x_1\right)\right)\right)\right) \\ &= V_p\left(x, J_p^B x_n - \lambda_n\left(Tx_n + \theta_n\left(J_p^B x_n - J_p^B x_1\right)\right)\right) \text{(by (21))} \\ &\leq V_p\left(x, J_p^B x_n\right) - p\lambda_n\left\langle J_p^{B^*}\left(J_p^B x_n - \lambda_n\left(Tx_n + \theta_n\left(J_p^B x_n - J_p^B x_1\right)\right)\right)\right. \\ &\quad \left. - x, Tx_n + \theta_n\left(J_p^B x_n - J_p^B x_1\right)\right\rangle \\ &= \phi_p(x, x_n) - p\lambda_n\left\langle x_n - x, Tx_n + \theta_n\left(J_p^B x_n - J_p^B x_1\right)\right\rangle \\ &\quad - p\lambda_n\left\langle J_p^{B^*}\left(J_p^B x_n - \lambda_n\left(Tx_n + \theta_n\left(J_p^B x_n - J_p^B x_1\right)\right)\right)\right. \\ &\quad \left. - x_n, Tx_n + \theta_n\left(J_p^B x_n - J_p^B x_1\right)\right\rangle.\end{aligned} \quad (33)$$

By Schwartz inequality and by applying inequality (32), we obtain

$$\begin{aligned}\phi_p(x, x_{n+1}) &\leq \phi_p(x, x_n) - p\lambda_n\left\langle x_n - x, Tx_n + \theta_n\left(J_p^B x_n - J_p^B x_1\right)\right\rangle \\ &\quad + p\lambda_n M_0 M \leq \phi_p(x, x_n) - p\lambda_n\langle x_n - x, Tx_n \\ &\quad - Tx\rangle \text{ (since } x \in T^{-1}(0)\text{)} \\ &\quad - p\lambda_n \theta_n\left\langle x_n - x, J_p^B x_n - J_p^B x_1\right\rangle + p\lambda_n M_0 M.\end{aligned} \quad (34)$$

By Lemma 6, $p\langle x - x_n, J_p^B x_n - J_p^B x_1\rangle \leq \phi_p(x, x_n) - \phi_p(x, x_1)$. Consequently, $p\langle x - x_n, J_p^B x_n - J_p^B x_1\rangle \leq \phi_p(x, x_n)$. Therefore, using $(p, \eta)$-strongly monotonicity property of $T$, we have

$$\begin{aligned}\phi_p(x, x_{n+1}) &\leq \phi_p(x, x_n) - p\eta\lambda_n\|x_n - x\|^p \\ &\quad - p\lambda_n\theta_n\left\langle x_n - x, J_p^B x_n - J_p^B x_1\right\rangle + p\lambda_n M_0 M \\ &\leq \phi_p(x, x_n) - p\lambda_n\|x_n - x\|^p \\ &\quad + p\lambda_n\theta_n\left\langle x - x_n, J_p^B x_n - J_p^B x_1\right\rangle + p\lambda_n M_0 M \\ &\leq \phi_p(x, x_n) - p\lambda_n\left(\phi_p(x, x_n) - \frac{p}{q}\|x\|^q\right) \\ &\quad + p\lambda_n\theta_n\phi_p(x, x_n) + p\lambda_n M_0 M \\ &= (1 - p\lambda_n)\phi_p(x, x_n) + p\lambda_n\left(\frac{p}{q}\|x\|^q\right) \\ &\quad + p\lambda_n\theta_n\phi_p(x, x_n) + p\lambda_n M_0 M \\ &\leq (1 - p\lambda_n)r + p\lambda_n\frac{r}{4} + p\lambda_n\frac{r}{2} + \frac{p\lambda_n}{4}r \\ &= \left(1 - p\lambda_n + p\lambda_n\frac{1}{4} + p\lambda_n\frac{1}{2} + p\lambda_n\frac{1}{4}\right)r = r.\end{aligned} \quad (35)$$

Hence, $\phi_p(x, x_{n+1}) \leq r$. By induction, $\phi_p(x, x_n) \leq r \forall n \in \mathbb{N}$. Thus, from inequality (20), $\{x_n\}_{n=1}^{\infty}$ is bounded.

*Part 2*: we now show that $\{x_n\}_{n=1}^{\infty}$ converges strongly to a solution of $Tx = 0$. $(p, \eta)$-strongly monotone implies a monotone and the range of $(J_p^B + tT)$ is all of $B^*$ for all $t > 0$. By Kohsaka and Takahashi [8], since $B$ is a reflexive smooth strictly convex space, we obtain for every $t > 0$ and $x \in B$, there exists a unique $x_t \in B$ such that

$$J_p^B x = J_p^B x_t + tTx_t. \quad (36)$$

Define $J_t^T x := x_t$; in other words, define a single-valued mapping $J_t^T : B \longrightarrow B$ by $J_t^T = (J_p^B + tT)^{-1} J_p^B$. Such a $J_t^T$ is called the resolvent of $T$. Setting $t := 1/\theta_n$ and by the result of Aoyama et al. [34] and Reich [35], for some $x_1 \in B$, there exists in $B$ a unique

$$y_n := \left(J_p^B + \frac{1}{\theta_n}T\right)^{-1} J_p^B x_1 \quad (37)$$

with $y_n \longrightarrow x \in T^{-1}(0)$. Obviously, one can obtain that

$$Ty_n = \theta_n \left(J_p^B x_1 - J_p^B y_n\right), \tag{38}$$

and $\{y_n\}_{n=1}^{\infty}$ is known to be bounded. Also it can be obtained that

$$\theta_n \left(J_p^B y_n - J_p^B x_1\right) + Ty_n = 0. \tag{39}$$

From (39), we have that

$$\theta_n \left(J_p^B y_n - J_p^B x_1\right) + Ty_n - Tx_1 = -Tx_1, \tag{40}$$

which is equivalent to

$$Ty_n - Tx_1 = -\theta_n \left(J_p^B y_n - J_p^B x_1\right) - Tx_1. \tag{41}$$

Consequently,

$$\begin{aligned}
\eta \|y_n - x_1\|^p &\le \langle y_n - x_1, Ty_n - Tx_1 \rangle \text{ (by } (p,\eta)\text{-strongly monotonicity of } T) \\
&= -\theta_n \langle y_n - x_1, \left(J_p^B y_n - J_p^B x_1\right) \rangle - \langle y_n - x_1, Tx_1 \rangle \\
&\le \|Tx_1\| \|y_n - x_1\|,
\end{aligned} \tag{42}$$

which shows that the sequence $\{y_n\}_{n=1}^{\infty}$ is bounded. Moreover, $\{x_n\}_{n=1}^{\infty}$ is bounded, and hence, $\{Tx_n\}_{n=1}^{\infty}$ is bounded. Following the same arguments as in part 1, we get

$$\begin{aligned}
\phi_p(y_n, x_{n+1}) &\le \phi_p(y_n, x_n) - p\lambda_n \langle x_n - y_n, Tx_n + \theta_n \left(J_p^B x_n - J_p^B x_1\right) \rangle \\
&\quad + p\lambda_n M_0 M.
\end{aligned} \tag{43}$$

By the $(p,\eta)$-strongly monotonicity property of $T$ and using Lemma 7 and Equation (39), we obtain

$$\begin{aligned}
&\langle x_n - y_n, Tx_n + \theta_n \left(J_p^B x_n - J_p^B x_1\right) \rangle \\
&= \langle x_n - y_n, Tx_n + \theta_n \left(J_p^B x_n - J_p^B y_n + J_p^B y_n - J_p^B x_1\right) \rangle \\
&= \theta_n \langle x_n - y_n, J_p^B x_n - J_p^B y_n \rangle \\
&\quad + \langle x_n - y_n, Tx_n + \theta_n \left(J_p^B y_n - J_p^B x_1\right) \rangle \\
&= \theta_n \langle x_n - y_n, J_p^B x_n - J_p^B y_n \rangle + \langle x_n - y_n, Tx_n - Ty_n \rangle \\
&\ge \theta_n g(\|x_n - y_n\|) + \eta \|x_n - y_n\|^p \ge \frac{1}{p} \theta_n \phi_p(y_n, x_n).
\end{aligned} \tag{44}$$

Therefore, the inequality (43) becomes

$$\phi_p(y_n, x_{n+1}) \le (1 - \lambda_n \theta_n) \phi_p(y_n, x_n) + p\lambda_n M_0 M. \tag{45}$$

Observe that by Lemma 6, we have

$$\begin{aligned}
\phi_p(y_n, x_n) &\le \phi_p(y_{n-1}, x_n) - p\langle y_n - x_n, J_p^B y_{n-1} - J_p^B y_n \rangle \\
&= \phi_p(y_{n-1}, x_n) + p\langle x_n - y_n, J_p^B y_{n-1} - J_p^B y_n \rangle \\
&\le \phi_p(y_{n-1}, x_n) + \|J_p^B y_{n-1} - J_p^B y_n\| \|x_n - y_n\|.
\end{aligned} \tag{46}$$

Let $R > 0$ such that $\|x_1\| \le R, \|y_n\| \le R$ for all $n \in \mathbb{N}$. We obtain from Equation (39) that

$$J_p^B y_{n-1} - J_p^B y_n + \frac{1}{\theta_n}(Ty_{n-1} - Ty_n) = \frac{\theta_{n-1} - \theta_n}{\theta_n}\left(J_p^B x_1 - J_p^B y_{n-1}\right). \tag{47}$$

By taking the duality pairing of each side of this equation with respect to $y_{n-1} - y_n$ and by the strong monotonicity of $T$, we have

$$\langle J_p^B y_{n-1} - J_p^B y_n, y_{n-1} - y_n \rangle \le \frac{\theta_{n-1} - \theta_n}{\theta_n} \|J_p^B x_1 - J_p^B y_{n-1}\| \|y_{n-1} - y_n\|. \tag{48}$$

Since $\{\theta_n\}_{n=1}^{\infty}$ is a decreasing sequence, it is known that $\theta_{n-1} \ge \theta_n$. Therefore,

$$\frac{\theta_{n-1} - \theta_n}{\theta_n} = \frac{\theta_{n-1}}{\theta_n} - 1 \ge 0. \tag{49}$$

Consequently,

$$\|J_p^B y_{n-1} - J_p^B y_n\| \le \left(\frac{\theta_{n-1}}{\theta_n} - 1\right) \|J_p^B y_{n-1} - J_p^B x_1\|. \tag{50}$$

Using (46) and (50), the inequality (45) becomes

$$\begin{aligned}
\phi_p(y_n, x_{n+1}) &\le (1 - \lambda_n \theta_n)\phi_p(y_{n-1}, x_n) + C\left(\frac{\theta_{n-1}}{\theta_n} - 1\right) \\
&\quad + p\lambda_n M_0 M,
\end{aligned} \tag{51}$$

for some constant $C > 0$. By Lemma 8, $\phi_p(y_{n-1}, x_n) \longrightarrow 0$ as $n \longrightarrow \infty$ and using Lemma 9, we have that $x_n - y_{n-1} \longrightarrow 0$ as $n \longrightarrow \infty$. Since $y_n \longrightarrow x \in T^{-1}(0)$, we obtain that $x_n \longrightarrow x$ as $n \longrightarrow \infty$.

**Corollary 14.** *Let $H$ be a Hilbert space, $p > 1, \eta \in (1,\infty)$ and suppose $T: H \longrightarrow H$ is a continuous, $(p,\eta)$-strongly monotone mapping such that $D(T) \subseteq \text{range}(I + tT)$ for all $t > 0$. For arbitrary $x_1 \in H$, define the sequence $\{x_n\}_{n=1}^{\infty}$ iteratively by*

$$x_{n+1} := x_n - \lambda_n Tx_n - \lambda_n \theta_n (x_n - x_1), \quad n \in \mathbb{N}, \tag{52}$$

*where $\{\lambda_n\}_{n=1}^{\infty} \subset (0,1)$ and $\{\theta_n\}_{n=1}^{\infty}$ in $(0,1/2)$ are real sequences satisfying the conditions:*

*(i) $\lim_{n\to\infty} \theta_n = 0$ and $\{\theta_n\}_{n=1}^{\infty}$ is decreasing*

(ii) $\sum_{n=1}^{\infty} \lambda_n \theta_n = \infty$

(iii) $\lim_{n\to\infty}((\theta_{n-1}/\theta_n) - 1)/\lambda_n\theta_n = 0$, $\sum_{n=1}^{\infty} \lambda_n < \infty \forall n \in \mathbb{N}$

Suppose that the equation $Tx = 0$ has a solution. Then, the sequence $\{x_n\}_{n=1}^{\infty}$ converges strongly to the solution of the equation $Tx = 0$.

*Proof.* The result follows from Theorem 13 since uniformly smooth and uniformly convex spaces are more general than the Hilbert spaces.

Examples are given for nonlinear mappings of the monotone type which satisfies the conditions stated in the main theorem.

*Example 15.* Let $T : \mathbb{R}^n \to \mathbb{R}^n$ with $Tx = |x|^{p-2}$ and $p \geq 2$. Then,

$$\langle x - y, |x|^{p-2}x - |y|^{p-2}y \rangle$$
$$= \frac{1}{2}(|x|^{p-2} + |y|^{p-2})|x - y|^2 + \frac{1}{2}(|x|^{p-2} - |y|^{p-2})(|x|^2 - |y|^2)$$
$$\geq 2^{-1}(|x|^{p-2} + |y|^{p-2})|x - y|^2 \geq 2^{2-p}|x - y|^p. \quad (53)$$

Thus, $T$ is $(p, \eta)$-strongly monotone with $\eta := 2^{2-p}$.

*Example 16.* Let $E = \mathbb{R}^2$ with the usual norm. Consider the function $T : \mathbb{R}^2 \to \mathbb{R}^2$ defined by $Tx = gx$, where

$$g = \begin{pmatrix} 8 & -5 \\ 5 & 13 \end{pmatrix},$$
$$x = \begin{pmatrix} x_1 \\ x_2 \end{pmatrix}. \quad (54)$$

Then, $T$ is $(p, \eta)$-strongly monotone with $p := 2$ and $\eta := 8$. Indeed, $\langle x, Tx \rangle \geq 8\|x\|^2$.

## 4. Solution of Convex Minimization Problems

The result of Theorem 13 is applied in this section for solving a problem of finding a minimizer of a convex function $\varphi$ defined from a real Banach space $B$ to $\mathbb{R}$. Recall that a mapping $T : B \to B^*$ is said to be *coercive* if for any $x \in B$,

$$\frac{\langle x, Tx \rangle}{\|x\|} \to \infty \quad \text{as} \|x\| \to \infty. \quad (55)$$

The following well-known basic results will be used.

**Lemma 17.** *Let $\varphi : B \to \mathbb{R}$ be a real-valued differentiable convex function and $u \in B$. Let $d\varphi : B \to B^*$ denote the differential map associated to $\varphi$. Then, the following hold:*

(1) *The point $u$ is a minimizer of $\varphi$ on $B$ if and only if $d\varphi(u) = 0$*

(2) *If $\varphi$ is bounded, then $\varphi$ is locally Lipschitzian, i.e., for every $x_0 \in B$ and $r > 0$, there exists $L > 0$ such that $\varphi$ is $L$-Lipschitzian on $B(x_0, r)$, i.e.,*

$$|\varphi(x) - \varphi(y)| \leq L\|x - y\| \quad \forall x, y \in B(x_0, r). \quad (56)$$

The main result in this section is given below.

**Theorem 18.** *Let $B$ be a uniformly smooth and uniformly convex real Banach space. Let $\varphi : B \to \mathbb{R}$ be a differentiable, convex, bounded, and coercive function. Let $\{\lambda_n\}_{n=1}^{\infty} \subset (0, 1)$ and $\{\theta_n\}_{n=1}^{\infty}$ in $(0, 1/2)$ be real sequences such that,*

(i) $\lim_{n\to\infty}\theta_n = 0$ and $\{\theta_n\}_{n=1}^{\infty}$ is decreasing

(ii) $\sum_{n=1}^{\infty} \lambda_n \theta_n = \infty$

(iii) $\lim_{n\to\infty}((\theta_{n-1}/\theta_n) - 1)/\lambda_n\theta_n = 0$, $\sum_{n=1}^{\infty} \lambda_n < \infty \forall n \in \mathbb{N}$

*For arbitrary $x_1 \in B$, define $\{x_n\}_{n=1}^{\infty}$ iteratively by*

$$x_{n+1} = J_p^{B^*}\left(J_p^B x_n - \lambda_n\left(d\varphi(x_n) + \theta_n\left(J_p^B x_n - J_p^B x_1\right)\right)\right), \quad n \in \mathbb{N}, \quad (57)$$

*where $J^B$ is the generalized duality mapping from $B$ into $B^*$. Then, $\varphi$ has a minimizer $x^* \in B$ and the sequence $\{x_n\}_{n=1}^{\infty}$ converges strongly to $x^*$.*

*Proof.* $\varphi$ has a minimizer because it is a function which is lower semicontinuous, convex, and coercive. Moreover, $x^* \in B$ minimizes $\varphi$ if and only if $d\varphi(x^*) = 0$. It can be inferred that $d\varphi$ is a maximal monotone due to the convexity, the differentiability, and the boundedness of $\varphi$ (see, e.g., Minty [36] and Moreau [37]). The next task is to show that $d\varphi$ is bounded. Indeed, let $x_0 \in B$ and $r > 0$. By Lemma 17, there exists $L > 0$ such that

$$|\varphi(x) - \varphi(y)| \leq L\|x - y\| \quad \forall x, y \in B(x_0, r). \quad (58)$$

Let $v^* \in d\varphi(B(x_0, r))$ and $x^* \in B(x_0, r)$ such that $v^* = d\varphi(x^*)$. Since $B(x_0, r)$ is open, for all $u \in B$, there exists $\sigma > 0$ such that $x^* + \sigma u \in (B(x_0, r)$. From the fact that $v^* = d\varphi(x^*)$ and inequality (58), it is obtained that

$$\langle v^*, \sigma u \rangle \leq \varphi(x^* + \sigma u) - \varphi(x^*) \leq \sigma L\|u\|, \quad (59)$$

such that

$$\langle v^*, u \rangle \leq L\|u\| \quad \forall u \in B. \quad (60)$$

Consequently, $\|v^*\| \leq L$, which implies that $d\varphi(B(x_0, r))$ is bounded. Thus, $d\varphi$ is bounded. Hence, it can be deduced from Theorem 13 that the sequence $\{x_n\}_{n=1}^{\infty}$ converges strongly to $x^*$, a minimizer of $\varphi$.

*Example 19.* An example of a function which is coercive is a real valued function $f : \mathbb{R}^2 \longrightarrow \mathbb{R}$ which is defined by $f(u, v) = u^4 - 7uv + v^3$.

Constructively, $f(u, v) = (u^3 + v^4)(1 - (7uv/(u^3 + v^4)))$. As $\|(u, v)\| \longrightarrow \infty$, $7uv/(u^3 + v^4) \longrightarrow 0$ while $u^3 + v^4 \longrightarrow \infty$. It follows that

$$\lim_{\|(u,v)\| \longrightarrow \infty} f(u, v) = \lim_{\|(u,v)\| \longrightarrow \infty} (u^3 + v^4)(1 - 0) = +\infty. \quad (61)$$

Hence, $f$ is coercive.

## 5. Solutions of Variational Inequality Problems

Let $K$ be a nonempty, closed, and convex subset of a real normed linear space $B$ and let $T : K \to B$ be a nonlinear mapping. The variational inequality problem is to

$$\text{find } x \in K \text{ such that } \langle j_p(x - y), Tx \rangle \geq 0, \quad \forall y \in K, \quad (62)$$

for some $j_p(x - y) \in J_p(x - y)$. The set of solutions of a variational inequality problem is denoted by $VI(T, K)$. If $B \coloneqq H$, a Hilbert space, the variational inequality problem reduces to

$$\text{find } x \in K \text{ such that } \langle x - y, Tx \rangle \geq 0, \quad \forall y \in K, \quad (63)$$

which was introduced and studied by Stampacchia [38]. Variational inequality theory has emerged as an important tool in studying a wide class of related problems arising in mathematical, physical, regional, engineering, and nonlinear optimization sciences. The theories of variational inequality problems have numerous applications in the study of nonlinear analysis (see, e.g., Censor et al. [39], Korpelevich [40], Shi [41], and Stampacchia [38] and the references contained in them). Several existence results have been established for (62) and (63) when $T$ is a monotone type mapping (see, e.g., Barbu and Precupanu [42], Browder [43], and Hartman and Stampacchia [44] and the references contained in them).

Let $K$ be a closed convex subset of $H$. The projection into $K$ is defined to be the mapping, $P_K : H \to K$, which is given by

$$\|P_K(x) - x\| = \min \{\|y - x\| : y \in K\}. \quad (64)$$

Gradient projection method is an orthodox way for solving (63). The projection algorithm is given by

$$\begin{cases} x_1 \in K, \\ x_{n+1} = P_K(x_n - \eta_n T(x_n)), \quad n \in \mathbb{N}, \end{cases} \quad (65)$$

where $T$ is $\eta$-strongly pseudomonotone and $L$-Lipschitz continuous mapping (see, e.g., Khanh and Vuong [45]). A recent report eliminated some drawbacks in the study of algorithm (65) [46]. The report considered a mapping $T$, which is $\eta$-strongly pseudomonotone and bounded on bounded subsets of $K$.

We are interested in the set of solutions of the form $VI(T, C)$, where $T : B \longrightarrow B^*$ is a $(p, \eta)$-strongly monotone mapping, $C \coloneqq \cap_{i=1}^{N} F(\varphi_i) \neq \varnothing, \varphi_i : K \longrightarrow B, i = 1, 2, \cdots, N$ is a finite family of quasi-$\phi_p$-nonexpansive mappings, and $B$ is a uniformly smooth and uniformly convex real Banach space. Recall that a mapping $\varphi : K \longrightarrow K$ is called nonexpansive if $\|\varphi x - \varphi y\| \leq \|x - y\|, \forall x, y \in K$. The set of fixed points of the mapping $\varphi$ will be denoted by $F(\varphi)$. A mapping $\varphi$ is said to be quasi-$\phi_p$-nonexpansive if $F(\varphi) \neq \varnothing$ and $\phi_p(u, \varphi x) \leq \phi_p(u, x), \forall x \in K$ and $u \in F(\varphi)$. The proof of the following theorem is given.

**Theorem 20.** *Let $B$ be a uniformly smooth and uniformly convex real Banach space and $K$ a nonempty, closed, and convex subset of $B$. Let $p > 1, \eta \in (1, \infty)$, suppose $T : B \longrightarrow B^*$ is a continuous, $(p, \eta)$-strongly monotone mapping such that the range of $(J_p + tT)$ is all of $B^*$ for all $t > 0$. Let $\varphi_i : K \longrightarrow B$, $i = 1, 2, \cdots, N$ be a finite family of quasi-$\phi_p$-nonexpansive mappings with $C \coloneqq \cap_{i=1}^{N} F(\varphi_i) \neq \varnothing$. Let $\{\lambda_n\}_{n=1}^{\infty} \subset (0, 1)$ and $\{\theta_n\}_{n=1}^{\infty}$ in $(0, 1/2)$ be real sequences such that*

(i) $\lim_{n \longrightarrow \infty} \theta_n = 0$ and $\{\theta_n\}_{n=1}^{\infty}$ is decreasing

(ii) $\sum_{n=1}^{\infty} \lambda_n \theta_n = \infty$

(iii) $\lim_{n \longrightarrow \infty} ((\theta_{n-1}/\theta_n) - 1)/\lambda_n \theta_n = 0, \sum_{n=1}^{\infty} \lambda_n < \infty \forall n \in \mathbb{N}$

*For arbitrary $x_1 \in B$, define $\{x_n\}_{n=1}^{\infty}$ iteratively by*

$$\begin{aligned} x_{n+1} = J_p^{B^*} \Big( J_p^B \big( \varphi_{[n]} x_n \big) - \lambda_n \Big( T \big( \varphi_{[n]} x_n \big) \\ + \theta_n \big( J_p^B \big( \varphi_{[n]} x_n \big) - J_p^B \big( \varphi_{[n]} x_1 \big) \big) \Big) \Big), \\ n \in \mathbb{N}, \end{aligned} \quad (66)$$

*where $\varphi_{[n]} \coloneqq \varphi_n \operatorname{Mod} N$ and $J_p^B$ is the generalized duality mapping from $B$ into $B^*$. Then, the sequence $\{x_n\}_{n=1}^{\infty}$ converges strongly to $x \in VI(T, C)$.*

*Proof.* Firstly, it is shown that the sequence $\{x_n\}_{n=1}^{\infty}$ is bounded.

Let $q > 1$ with $1/p + 1/q = 1$ and $x \in VI(T, C)$. It suffices to show that $\phi_p(x, x_n) \leq r, \forall n \in \mathbb{N}$. The proof is by induction. Let $r > 0$ be sufficiently large such that

$$r \geq \max \left\{ \phi_p(x, x_1), 4M_0 M, \frac{4p}{q} \|x\|^q \right\}, \quad (67)$$

where $M_0 > 0$ and $M > 0$ are arbitrary but fixed. By construction, $\phi_p(x, x_1) \leq r$. Suppose that $\phi_p(x, x_n) \leq r$ for some $n \in \mathbb{N}$. From inequality (20), for real $p > 1$, we have $\|x_n\| \leq r^{1/p} + \|x\|$. Let $D \coloneqq \{z \in B : \phi_p(x, z) \leq r\}$. We show that $\phi_p(x, x_{n+1}) \leq r$. It is known that $T$ is locally bounded and $J_p^B$ is uniformly continuous on bounded subsets of $B$.

Define

$$M_0 := \sup\left\{\|T(\varphi_{[n]}x_n) + \theta_n(J_p^B(\varphi_{[n]}x_n) - J_p^B(\varphi_{[n]}x_1))\| : \theta_n \in \left(0, \frac{1}{2}\right), x_n \in D\right\} + 1. \quad (68)$$

Let $\psi$ denotes the modulus of continuity of $J_p^{B^*}$. Then,

$$\begin{aligned}\|\varphi_{[n]}x_n - x_{n+1}\| &= \|\varphi_{[n]}x_n - J_p^{B^*}(J_p^B(\varphi_{[n]}x_n) - \lambda_n(T(\varphi_{[n]}x_n) \\ &\quad + \theta_n(J_p^B(\varphi_{[n]}x_n) - J_p^B(\varphi_{[n]}x_1))))\| \\ &= \|J_p^{B^*}(J_p^B(\varphi_{[n]}x_n)) - J_p^{B^*}(J_p^B(\varphi_{[n]}x_n) \\ &\quad - \lambda_n(T(\varphi_{[n]}x_n) + \theta_n(J_p^B(\varphi_{[n]}x_n) \\ &\quad - J_p^B(\varphi_{[n]}x_1))))\| \le \psi(|\lambda_n|\|T(\varphi_{[n]}x_n) \\ &\quad + \theta_n(J_p^B(\varphi_{[n]}x_n) - J_p^B(\varphi_{[n]}x_1))\|) \\ &\le \psi(|\lambda_n|M_0) \le \psi(\sup\{|\lambda_n|M_0 : \lambda_n \in (0,1)\}).\end{aligned} \quad (69)$$

Since $T$ is locally bounded and the duality mapping $J_p^B$ is uniformly continuous on bounded subsets of $B$, the $\sup\{|\lambda_n|M_0\}$ exists, and it is a real number different from infinity. Define $M := \psi(\sup\{|\lambda_n|M_0\})$. Applying Lemma 4 with $y^* := \lambda_n(T(\varphi_{[n]}x_n) + \theta_n(J_p^B(\varphi_{[n]}x_n) - J_p^B(\varphi_{[n]}x_1)))$ and by using the definition of $x_{n+1}$, we compute as follows:

$$\begin{aligned}\phi_p(x, x_{n+1}) &= \phi_p(x, J_p^{B^*}(J_p^B(\varphi_{[n]}x_n) - \lambda_n(T(\varphi_{[n]}x_n) \\ &\quad + \theta_n(J_p^B(\varphi_{[n]}x_n) - J_p^B(\varphi_{[n]}x_1))))) \\ &= V_p(x, J_p^B(\varphi_{[n]}x_n) - \lambda_n(T(\varphi_{[n]}x_n) \\ &\quad + \theta_n(J_p^B(\varphi_{[n]}x_n) - J_p^B(\varphi_{[n]}x_1)))) \text{ (by (21))} \\ &\le V_p(x, J_p^B(\varphi_{[n]}x_n)) - p\lambda_n\langle x_{n+1} - x, T(\varphi_{[n]}x_n) \\ &\quad + \theta_n(J_p^B(\varphi_{[n]}x_n) - J_p^B(\varphi_{[n]}x_1))\rangle \\ &= \phi_p(x, \varphi_{[n]}x_n) - p\lambda_n\langle \varphi_{[n]}x_n - x, T(\varphi_{[n]}x_n) \\ &\quad + \theta_n(J_p^B(\varphi_{[n]}x_n) - J_p^B(\varphi_{[n]}x_1))\rangle \\ &\quad - p\lambda_n\langle x_{n+1} - \varphi_{[n]}x_n, T(\varphi_{[n]}x_n) \\ &\quad + \theta_n(J_p^B(\varphi_{[n]}x_n) - J_p^B(\varphi_{[n]}x_1))\rangle.\end{aligned} \quad (70)$$

By Schwartz inequality and by applying inequality (69), we obtain

$$\begin{aligned}\phi_p(x, x_{n+1}) &\le \phi_p(x, \varphi_{[n]}x_n) - p\lambda_n\langle \varphi_{[n]}x_n - x, T(\varphi_{[n]}x_n) \\ &\quad + \theta_n(J_p^B(\varphi_{[n]}x_n) - J_p^B(\varphi_{[n]}x_1))\rangle + p\lambda_nM_0M \\ &\le \phi_p(x, \varphi_{[n]}x_n) - p\lambda_n\langle \varphi_{[n]}x_n - x, T(\varphi_{[n]}x_n) - Tx\rangle \\ &\quad - p\lambda_n\langle \varphi_{[n]}x_n - x, Tx\rangle - p\lambda_n\theta_n\langle \varphi_{[n]}x_n - x, J_p^B(\varphi_{[n]}x_n) \\ &\quad - J_p^B(\varphi_{[n]}x_1)\rangle + p\lambda_nM_0M \le \phi_p(x, \varphi_{[n]}x_n) \\ &\quad - p\lambda_n\langle \varphi_{[n]}x_n - x, T(\varphi_{[n]}x_n) - Tx\rangle (\text{since } x \in VI(T, C)) \\ &\quad - p\lambda_n\theta_n\langle \varphi_{[n]}x_n - x, J_p^B(\varphi_{[n]}x_n) - J_p^B(\varphi_{[n]}x_1)\rangle \\ &\quad + p\lambda_nM_0M.\end{aligned} \quad (71)$$

By Lemma 6, $p\langle x - \varphi_{[n]}x_n, J_p^B(\varphi_{[n]}x_n) - J_p^B(\varphi_{[n]}x_1)\rangle \le \phi_p(x, \varphi_{[n]}x_n) - \phi_p(x, \varphi_{[n]}x_1)$. Consequently, $p\langle x - \varphi_{[n]}x_n, J_p^B(\varphi_{[n]}x_n) - J_p^B(\varphi_{[n]}x_1)\rangle \le \phi_p(x, \varphi_{[n]}x_n)$. Therefore, using $(p, \eta)$-strongly monotonicity property of $T$, we have

$$\begin{aligned}\phi_p(x, x_{n+1}) &\le \phi_p(x, \varphi_{[n]}x_n) - p\eta\lambda_n\|T(\varphi_{[n]}x_n) - x\|^p \\ &\quad - p\lambda_n\theta_n\langle \varphi_{[n]}x_n - x, J_p^B(\varphi_{[n]}x_n) - J_p^B(\varphi_{[n]}x_1)\rangle \\ &\quad + p\lambda_nM_0M \le \phi_p(x, \varphi_{[n]}x_n) \\ &\quad - p\lambda_n\left(\phi_p(x, T(\varphi_{[n]}x_n)) - \frac{p}{q}\|x\|^q\right) \\ &\quad + p\lambda_n\theta_n\langle \varphi_{[n]}x_n - x, J_p^B(\varphi_{[n]}x_n) \\ &\quad - J_p^B(\varphi_{[n]}x_1)\rangle + p\lambda_nM_0M \\ &\le \phi_p(x, \varphi_{[n]}x_n) - p\lambda_n\phi_p(x, \varphi_{[n]}x_n) + p\lambda_n\frac{p}{q}\|x\|^q \\ &\quad + p\lambda_n\theta_n\phi_p(x, \varphi_{[n]}x_n) + p\lambda_nM_0M \\ &= (1 - p\lambda_n)\phi_p(x, \varphi_{[n]}x_n) + p\lambda_n\left(\frac{p}{q}\|x\|^q\right) \\ &\quad + p\lambda_n\theta_n\phi_p(x, \varphi_{[n]}x_n) + p\lambda_nM_0M \\ &\le (1 - p\lambda_n)\phi_p(x, x_n) + p\lambda_n\left(\frac{p}{q}\|x\|^q\right) \\ &\quad + p\lambda_n\theta_n\phi_p(x, x_n) + p\lambda_nM_0M \\ &\le (1 - p\lambda_n)r + p\lambda_n\frac{r}{4} + p\lambda_n\frac{r}{2} + \frac{p\lambda_n}{4}r \\ &= \left(1 - p\lambda_n + p\lambda_n\frac{1}{4} + p\lambda_n\frac{1}{2} + p\lambda_n\frac{1}{4}\right)r = r.\end{aligned} \quad (72)$$

Hence, $\phi_p(x, x_{n+1}) \le r$. By induction, $\phi_p(x, x_n) \le r \forall n \in \mathbb{N}$. Thus, from inequality (20), $\{x_n\}_{n=1}^\infty$ is bounded. The remaining part of the proof follows from the proof of Theorem 13.

*Remark 21.* It well known that uniformly smooth and uniformly convex spaces are more general than the Hilbert spaces. Therefore, the following corollary is readily obtainable.

**Corollary 22.** *Let $H$ be a Hilbert space and and $K$ a nonempty, closed, and convex subset of $H$. Let $p > 1, \eta \in (1,\infty)$; suppose $T : H \to H$ is a continuous, $(p,\eta)$-strongly monotone mapping such that $D(T) \subseteq \text{range}(I + tT)$ for all $t > 0$. Let $\varphi_i : K \to H$, $i = 1, 2, \cdots, N$ be a finite family of quasi-$\phi_p$-nonexpansive mappings with $C := \cap_{i=1}^{N} F(\varphi_i) \neq \emptyset$. Let $\{\lambda_n\}_{n=1}^{\infty} \subset (0, 1)$ and $\{\theta_n\}_{n=1}^{\infty}$ in $(0, 1/2)$ be real sequences such that*

(i) $\lim_{n\to\infty} \theta_n = 0$ and $\{\theta_n\}_{n=1}^{\infty}$ is decreasing

(ii) $\sum_{n=1}^{\infty} \lambda_n \theta_n = \infty$

(iii) $\lim_{n\to\infty}((\theta_{n-1}/\theta_n) - 1)/\lambda_n \theta_n = 0$, $\sum_{n=1}^{\infty} \lambda_n < \infty \forall n \in \mathbb{N}$

*For arbitrary $x_1 \in H$, define $\{x_n\}_{n=1}^{\infty}$ iteratively by*

$$x_{n+1} = \varphi_{[n]} x_n - \lambda_n \left( T\left(\varphi_{[n]} x_n\right) + \theta_n \left(\varphi_{[n]} x_n - \varphi_{[n]} x_1\right)\right), \quad n \in \mathbb{N}, \tag{73}$$

*where $\varphi_{[n]} := \varphi_n \text{Mod} N$. Then, the sequence $\{x_n\}_{n=1}^{\infty}$ converges strongly to $x \in VI(T, C)$.*

## 6. Conclusion

Real-life problems are usually modeled by nonlinear equations. Nonlinear equations occur in modeling problems, such as minimizing costs in industries and minimizing risks in businesses. Nonlinear equations of $(p, \eta)$-strongly monotone type, where $\eta \in (1,\infty), p > 1$, have been studied in this paper. The result was applied to obtain the solution of convex minimization and variational inequality problems, which have applications in several fields such as economics, game theory, and the sciences.

## Data Availability

Data sharing is not applicable to this article as no datasets were generated or analyzed during the current study.

## Disclosure

The abstract of this manuscript is submitted for presentation in the "9th International Conference on Mathematical Modeling in Physical Sciences," September 7–10, 2020, Tinos island, Greece.

## Conflicts of Interest

The authors declare no conflicts of interest.

Abstract and Applied Analysis                                                                                                     11problems," *Fixed Point Theory and Applications*, vol. 2016, no. 1, Article ID 97, 2016.

[19] M. O. Aibinu and O. T. Mewomo, "Algorithm for the generalized $\Phi$-strongly monotone mappings and application to the generalized convex optimization problem," *Proyecciones (Antofagasta)*, vol. 38, no. 1, pp. 59–82, 2019.

[20] F. Altomare and M. Campiti, *Korovkin-Type Approximation Theory and Its Applications*, de Gruyter, Berlin: New York, 1994.

[21] O. Forster, *Analysis 1, differential- und Integralrechnung einer veränderlichen*, Vieweg Studium, 1983.

[22] F. Kohsaka and W. Takahashi, "Existence and approximation of fixed points of firmly nonexpansive type mappings in Banach spaces," *SIAM Journal on Optimization*, vol. 19, no. 2, pp. 824–835, 2008.

[23] I. Cioranescu, *Geometry of Banach Spaces, Duality Mappings and Nonlinear Problems*, Kluwer Academic Publishers Group, Dordrecht, 1990.

[24] Z. B. Xu and G. F. Roach, "Characteristic inequalities of uniformly convex and uniformly smooth Banach spaces," *Journal of Mathematical Analysis and Applications*, vol. 157, no. 1, pp. 189–210, 1991.

[25] C. Zalinescu, "On uniformly convex functions," *Journal of Mathematical Analysis and Applications*, vol. 95, no. 2, pp. 344–374, 1983.

[26] Y. Alber, "Metric and generalized projection operators in Banach spaces: properties and applications," in *Theory and Applications of Nonlinear Operators of Accretive and Monotone Type*, A. G. Kartsatos, Ed., pp. 15–50, Marcel Dekker, New York, 1996.

[27] E. H. Zarantonello, "The meaning of the Cauchy-Schwarz-Buniakovsky inequality," *Proceedings of the American Mathematical Society*, vol. 59, no. 1, pp. 133–137, 1976.

[28] H. K. Xu, "Inequalities in Banach spaces with applications," *Nonlinear Analysis*, vol. 16, no. 12, pp. 1127–1138, 1991.

[29] H. K. Xu, "Iterative algorithms for nonlinear operators," *Journal of the London Mathematical Society*, vol. 66, no. 1, pp. 240–256, 2002.

[30] S. Kamimura and W. Takahashi, "Strong convergence of a proximal-type algorithm in a Banach space," *SIAM Journal on Optimization*, vol. 13, no. 3, pp. 938–945, 2002.

[31] R. T. Rockafellar, "Local boundedness of nonlinear, monotone operators," *Michigan Mathematical Journal*, vol. 16, no. 4, pp. 397–407, 1969.

[32] D. Pascali and S. Sburlan, *Nonlinear Mappings of Monotone Type*, Editura Academiae, Bucharest, Romania, 1978.

[33] C. E. Chidume and N. Djitte, "Strong convergence theorems for zeros of bounded maximal monotone nonlinear operators," *Abstract and Applied Analysis*, vol. 2012, Article ID 681348, 19 pages, 2012.

[34] K. Aoyama, F. Kohsaka, and W. Takahashi, "Proximal point methods for monotone operators in Banach Spaces," *Taiwanese Journal of Mathematics*, vol. 15, no. 1, pp. 259–281, 2011.

[35] S. Reich, "Constructive techniques for accretive and monotone operators," in *Applied Nonlinear Analysis*, pp. 335–345, Academic Press, New York, 1979.

[36] G. J. Minty, "Monotone (nonlinear) operators in Hilbert space," *Duke Mathematical Journal*, vol. 29, no. 3, pp. 341–346, 1962.

[37] J. J. Moreau, "Proximité et dualité dans un espace hilbertien," *Bulletin de la Société mathématique de France*, vol. 79, pp. 273–299, 1951.

[38] G. Stampacchia, "Formes bilineaires coercitives sur les ensembles convexes," *Comptes Rendus de l'Academie des Sciences*, vol. 258, pp. 4413–4416, 1964.

[39] Y. Censor, A. Gibali, and S. Reich, "Extensions of Korpelevich's extragradient method for the variational inequality problem in Euclidean space," *Optimization*, vol. 61, no. 9, pp. 1119–1132, 2012.

[40] G. M. Korpelevich, "An extragradient method for finding saddle points and for other problems," *Economics and Mathematical Methods*, vol. 12, pp. 747–756, 1967.

[41] P. Shi, "Equivalence of variational inequalities with Wiener–Hopf equations," *Proceedings of American Mathematical Society*, vol. 111, no. 2, pp. 339–346, 1991.

[42] V. Barbu and T. Precupanu, *Convexity and Optimization in Banach Spaces*, Springer, New York, 4th edition, 2012.

[43] F. E. Browder, "Nonlinear monotone operators and convex sets in Banach spaces," *Bulletin of the American Mathematical Society*, vol. 71, no. 5, pp. 780–786, 1965.

[44] P. Hartman and G. Stampacchia, "On some non-linear elliptic differential functional equations," *Acta Mathematica*, vol. 115, pp. 271–310, 1966.

[45] P. D. Khanh and P. T. Vuong, "Modified projection method for strongly pseudomonotone variational inequalities," *Journal of Global Optimization*, vol. 58, no. 2, pp. 341–350, 2014.

[46] T. N. Hai, "On gradient projection methods for strongly pseudomonotone variational inequalities without Lipschitz continuity," *Optimization Letters*, vol. 14, no. 5, pp. 1177–1191, 2020.